\documentclass[12pt,reqno, oneside]{amsart}
\makeatletter
\@namedef{subjclassname@2020}{%
  \textup{2020} Mathematics Subject Classification}
\makeatother
\usepackage{amsthm, mathrsfs, amsmath, amstext, esint, amsxtra, amsfonts, slashed, dsfont, amssymb, amsrefs}
\usepackage{verbatim,cite,relsize,color,soul}
\usepackage[latin1]{inputenc}
\usepackage{microtype}
\usepackage{color,enumitem,graphicx, xcolor}
\usepackage{mathtools}
\mathtoolsset{showonlyrefs}
\usepackage{bm} %bold greek letters
\usepackage[colorlinks=true,urlcolor=blue, citecolor=red,linkcolor=blue,
linktocpage,pdfpagelabels, bookmarksnumbered,bookmarksopen]{hyperref}
\usepackage[english]{babel}
\usepackage{geometry}

\numberwithin{equation}{section}
\newtheorem{theorem}{Theorem}[section]
\newtheorem{lemma}[theorem]{Lemma}

\newtheorem{corollary}[theorem]{Corollary}
\newtheorem{prop}[theorem]{Proposition}
\theoremstyle{definition}

\newtheorem{remark}[theorem]{Remark}

\newcommand{\R}{\mathbb R}
\newcommand{\C}{\mathbb C}

\newcommand{\eps}{\epsilon}
\DeclareMathOperator*{\loc}{loc}

\title[Subcritical Half-Wave Equation mixed nonlinearities]{Mass-subcritical Half-Wave Equation with mixed nonlinearities: existence and non-existence of ground states}
\author[J. Bellazzini, L. Forcella
]{Jacopo Bellazzini and Luigi Forcella}

\address[J. Bellazzini]{Department of Mathematics, University of Pisa, Largo Bruno Pontecorvo, 5, 56127, Pisa, Italy}
\email{jacopo.bellazzini@unipi.it}

\address[L. Forcella]{Department of Mathematics, University of Pisa, Largo Bruno Pontecorvo, 5, 56127, Pisa, Italy}
\email{luigi.forcella@unipi.it}



\begin{document}

%%%%%%%%%%%%%%%%%%%%%%%%%%%	

\begin{abstract}
We consider the problem of existence of constrained minimizers for the focusing mass-subcritical Half-Wave equation with a defocusing mass-subcritical perturbation. We show the existence of a critical mass such that minimizers do exist for any mass larger than or equal to the critical one, and do not exist below it. At the dynamical level, in the one dimensional case, we show that the ground states are orbitally stable.
\end{abstract}
	
%%%%%%%%%%%%%%%%%%%%%%%%%%%	
	
\maketitle

%%%%%%%%%%%%%%%%%%%%%%%%%

\section{Introduction} 
In this paper, we are motivated by  the following  Half-Wave type equation
\begin{equation}\label{SP}
i\psi_{t}= \sqrt{-\Delta} \psi +|\psi|^{q-1} \psi-|\psi|^{p-1}\psi,
\end{equation} 
where $\psi(t,x):[0,T)\times\R^{d}\rightarrow \C$, with $d\in\mathbb N$ and $T>0$, is a complex wave function and the exponents $q, p$ are of mass subcritical type, meaning that $1<q< p<1+\frac{2}{d}$. The operator $\sqrt{-\Delta}$ in \eqref{SP} is defined by the Fourier multiplier $m(\xi)=|\xi|$, $\xi\in\R^d$, and it acts as $\sqrt{-\Delta}f=\mathcal F^{-1}(|\xi|\mathcal F(f))$, where $\mathcal F$ and $\mathcal F^{-1}$ are the Fourier transform and its inverse, respectively. \medskip

In the 1D case, given an initial datum $\psi_0\in H^{1/2}(\R)$ associated to \eqref{SP}, it can be proved in the same manner as in \cite{KLR} (see also Section 3 for details) that the Cauchy problem is (globally) well-posed and satisfies the conservation of mass, energy, and momentum, defined  by 
\begin{equation}\label{eq:mass}
    \mathcal M(\psi(t))=\int_{\R}|\psi(t,x)|^2dx
\end{equation}

\begin{equation}\label{eq:en}
\begin{aligned}
    \mathcal E(\psi(t))&=\frac12\int_{\R} \bar{\psi}(t,x)\sqrt{-\partial_x^2} \psi(t,x)dx\\
&\phantom{=}    +\frac{1}{q+1}\int_{\R}| \psi(t,x)|^{q+1}dx-\frac{1}{p+1}\int_{\R}| \psi(t,x)|^{p+1}dx,
    \end{aligned}
\end{equation} 
and
\begin{equation}\label{eq:mom}
    \mathcal P(\psi(t))=-\int_{\R} i \bar \psi(t,x) \partial_x\psi(t,x)dx,
\end{equation} 
respectively. Particular solutions to \eqref{SP}
are the so-called \emph{traveling wave solutions}, which are  solutions to \eqref{SP} of the form
 \begin{equation}\label{eq:tv}
 \psi(t,x)=e^{i\omega t}u(x-vt)
 \end{equation}
 with $\omega\in \R$, $v
 \in\mathbb R$ such that $|v|< 1$, and $u(x)\in \C$ is a time-independent function belonging to  $H^{1/2}(\R)$ which satisfies 
\begin{equation}\label{eq:stat-1d}
\sqrt{- \partial_x^2} u +iv\partial_x u + \omega u +|u|^{q-1} u-|u|^{p-1}u=0.
\end{equation}

\noindent From a time-independent point of view, we can actually consider \eqref{eq:stat-1d} in arbitrary dimension, namely we consider solutions $u(x)\in \C$, $x\in\R^d$, to 
\begin{equation}\label{eq}
\sqrt{- \Delta} u +iv\cdot \nabla u + \omega u +|u|^{q-1} u-|u|^{p-1}u=0,
\end{equation}
where $\omega\in\R$ and $v\in\R^d$ with $|v|<1$. Motivated by the 1D case, since in that setting the mass is a physical quantity which is preserved along the flow of \eqref{SP}, a natural way  to find solutions $u$ to \eqref{eq}  is to look for  critical points of the  functional 
 \[
 E_v(u)=\frac12T_v(u)+\frac{1}{q+1}\int_{\R^{d}}| u|^{q+1}dx-\frac{1}{p+1}\int_{\R^{d}}|u|^{p+1}dx
 \]
 where 
 \[
T_v(u) =\int_{\R^d} \bar{u}(x)(\sqrt{-\Delta}+iv\cdot\nabla) u(x)dx,
 \]
constrained on the $L^{2}$-spheres of $H^{1/2}(\R^{d})$ described by 
\[S_{\rho}=\{u\in H^{1/2}(\R^{d}) \ \hbox{s.t.} \   \|u\|_{L^2}=\rho\}.
\]
 So, by a solution of \eqref{eq} we mean  a couple 
 $(\omega_{\rho}, u_\rho)\in\R\times H^{1/2}(\R^{d})$ where $\omega_{\rho}$ appears as the Lagrange multiplier associated to the critical point $u_{\rho}$ on $S_{\rho}$. Note that the functional $E_v(\psi)$ of the time-dependent function $\psi$ defined in \eqref{eq:tv} is preserved along the flow on \eqref{SP}
 by the fact that
 \[
 E_v(\psi)=\mathcal E(\psi)-\frac{v}{2}\mathcal {P}(\psi).\]

For a fixed a mass $\rho>0$, we introduce the \emph{ground state energy} as the quantity $I_{\rho^2}$, defined as
\begin{equation}\label{mini1}
I_{\rho^2}=\inf_{S_{\rho}} E_v(u).
\end{equation}

\noindent If a minimizer $u_\rho$ to \eqref{mini1} exists, we call it \emph{ground state solution}. 
Another functional which will play a relevant role for our approach is the following Pohozaev functional
\[ G_v(u)=T_v(u)+\frac{d(q-1)}{q+1}\int_{\R^{d}}| u|^{q+1}dx-\frac{d(p-1)}{p+1}\int_{\R^{d}}|u|^{p+1}dx.
\]
It is well-known that if $u_\rho$ is a solution to \eqref{eq}
then $u_\rho \in V_{\rho}$ 
where 
\[
V_{\rho}=\{u\in H^{1/2}(\R^{d})\ \hbox{s.t.} \  \|u\|_{L^2(\R^d)}=\rho, \ G_v(u)=0\}.
\]
The  goal of our paper is to establish a result about the existence and the non-existence of solutions to the minimization problem \eqref{mini1}. Specifically, for any fixed $v\in\R^d$ with $|v|<1$, the Theorem below shows the existence of a critical mass $\rho_0^v$ such that ground states exist for any mass $\rho\geq \rho_0^v$, and do not exist if $\rho<\rho_0^v$. \medskip

Our main result is as follows. 
\begin{theorem}\label{mainthm}
 Let $d\in\mathbb N$, $v\in \R^d$ with $|v|<1$, and $1<q< p<1+\frac{2}{d}$, then there exists a strictly positive mass $\rho_0^v$ such that:\\
 \noindent \textup{(i)} $ I_{\rho^2}=0$  for all  $\rho\in (0,\rho_0^v]$;\\
\noindent \textup{(ii)} $ I_{\rho^2}<0$  for all  $\rho\in (\rho_0^v, \infty)$.\\
\noindent Moreover, there are no constrained minimizers for $0<\rho<\rho_0^v$, and  for all  $\rho\in [\rho_0^v, \infty)$ there  exists $u_{\rho}\in S_{\rho}$ such that $I_{\rho^2}=E(u_{\rho})$. 
\end{theorem}

We give some comments about Theorem \ref{mainthm}.

\begin{remark}
\noindent \textup{(i)} The strategy to prove Theorem \ref{mainthm} follows a general scaling argument inspired by \cite{BS-JFA}, and, for $v=0$, our approach remains valid for any fractional NLS equation with mass-subcritical nonlinearities of power-type as appearing in \eqref{SP}. \smallskip 

\noindent \textup{(ii)} We also mention the paper by Jeanjean and Luo in \cite{JL-ZAMP} for the Schr\"odinger-Poisson equation, which was of inspiration for our paper concerning the non-existence results. We recall that in the Schr\"odinger-Poisson equation \[
-\Delta u+\omega u+(|x|^{-1}\ast |u|^2)u-|u|^{p-1}u=0,
\] 
the nonlocality is in the nonlinearity. On the other hand, the nonlocal nature of  equation  \eqref{SP} is due to the linear operator defining the kinetic energy.   In \cite{JL-ZAMP}, the analysis is specifically tailored to treat the nonlocal Coulomb-type term, while in our paper we are concerned with local-type nonlinearities in the whole mass-subcritical regime, and our achievements extend to any fractional Laplace operator.
\end{remark}
\begin{remark}
We emphasize that when $\rho=\rho_0^v$ the weak subadditivity inequality $I_{(\rho_0^v)^2}\leq I_{\mu^2}+I_{(\rho_0^v)^2-\mu^2}$ that always holds for this class of translation invariant minimization problem, becomes an equality
by the fact that for $0<\mu<\rho_0^v$ we have $I_{\mu^2}=I_{(\rho_0^v)^2}=0$. This fact shows that in general the strong subadditivity inequality is only a sufficient condition for the existence of constrained minimizers.
\end{remark}

\begin{remark} \noindent \textup{(i)} After the seminal works  by Tao, Vi\c san, and Zhang  \cite{TVZ-CPDE} on the classical NLS equation with combined nonlinearities, the interest on this type of equation rapidly increased. 
In particular, after the  papers by Soave \cite{So-JDE, So-JFA}, several scholars  treated the problem of the  existence of ground states with fixed mass for NLS-type equations with mixed nonlinearities. We specifically mention the paper by the authors   \cite{BFG-Annali},  and by  Jeanjean and Lu \cite{JL-2020}, that completed the study on the different scenarios initiated by Soave. We also cite the paper by Jeanjean and Lu \cite{JL-2022} that treated general nonlinearities under suitable growth conditions. \\
We remark that following verbatim the approach in our paper for the classical NLS equation with two competing nonlinearities,   we can recover the result of \cite{JL-2022}. It is worth mentioning, however, that the authors of \cite{JL-2022} consider general mass-subcritical nonlinearities and not only power-like terms.\smallskip 

\noindent \textup{(ii)} Concerning the Half-Wave equation with two competing nonlinearities, we should mention the paper \cite{Z-L}, where the authors study the problem of existence and  non-existence of traveling waves under the mass constraint. The defocusing-focusing case in the mass-subcritical case (i.e., under the same assumption of Theorem \ref{mainthm}) is considered in  \cite[Theorem 1.2]{Z-L}. The authors of that paper claim that ground states cannot exist with sufficiently small mass, as well as they claim the non-existence of a ground state with zero energy. However, the  proof of  \cite[Theorem 1.2]{Z-L} contains a crucial mistake in the sign of some estimates when treating  the defocusing term. Furthermore, we actually prove in Theorem \ref{mainthm} that a ground state with zero energy do exists, and its corresponding mass is exactly the threshold for the existence of ground states.\smallskip

\noindent \textup{(iii)} As for the existence of ground states for \eqref{SP} in the mass supercritical regime, we refer to \cite{ZLZ-NA}.

\end{remark}

\begin{remark}
Concerning  the existence of ground states, the case with only one focusing term is quite easy, due to the scaling invariance of the equation.  On the other hand, the problem becomes very interesting and more difficult  from the perspective of the uniqueness and the non-degeneracy of the ground state when $v=0$, see 
\cite{FL-acta, FLS-cpam, AT-acta}, and also regarding the lack of small data scattering on the dynamical side when $v\neq 0$, see \cite{BGLV}.
\end{remark}

Our second result concerns the dynamics of the ground states, and specifically it is the following orbital stability property of standing and traveling waves in 1D. 
\begin{corollary}\label{thm:corol} Fix $d=1$,  $v\in \R$ with $|v|<1$, and  let $\rho_0^v$ be as given in Theorem \ref{mainthm}. For any $\rho>\rho_0^v$, the set 
\[
\mathcal G=\{e^{i\gamma}u(\cdot+y) \ \hbox{s.t.} \   \gamma\in\R,  \ y\in \R, \  \hbox{ and } \ u\in S_\rho \ \hbox{ with } \ I_{\rho^2}=E_v(u)\}
\]
 is orbitally stable.
\end{corollary}
\begin{remark} It is worth noticing that at a dynamical level we can only consider the 1D framework, as is the only case where we can state a local well-posedness theory of the Cauchy problem associated to \eqref{SP}. 
\end{remark}

\subsection{Notations} We will work in the whole domain $\R^d$, hence we will omit the dependence on it, and systematically write norms, integrals, etc. without mention of the space we are working on. We will use the short notation $\|f\|_p$ for the $L^p$-norm of a function $f$. $H^s$, for $s\in(0,1)$, is the standard Sobolev space endowed with norm $\|f\|_{H^s}=\|(1-\Delta)^{s/2}f\|_2$. $\Re z$, $\Im z$, and $\bar z$ are the real part, the imaginary part, and the complex conjugate of $z$, respectively. If there is no confusion, we suppress the subscript from the functionals $E_v,G_v,T_v$ and the superscript from $
\rho_0^v$ from now on.

%%%%%%%%%%%%%%%%%%%%%%%%

\section{Proof of the main result} 
 
As described in the Introduction, the aim is to study the existence of minimizers for  
 \[
 E(u)=\frac12T(u)+\frac{1}{q+1}\int | u|^{q+1}dx-\frac{1}{p+1}\int |u|^{p+1}dx
 \]
 where 
 \[
T(u) =\int \bar{u}(x)(\sqrt{-\Delta}+iv\cdot\nabla) u(x)dx,
 \]
for $|v|<1$, under  mass constraint  
\[S_{\rho}=\{u\in H^{1/2} \ \hbox{s.t.} \   \|u\|_{2}=\rho\},
\]
for $1<q<p<1+\frac2d$.
Then, the problem is to compute
\[
I_{\rho^2}=\inf_{S_{\rho}}E(u),
\]
which is equivalent to 
\[
I_{\rho^2}=\inf_{V_{\rho}}E(u),
\]
where $V_\rho=\{ u\in H^1 \ \hbox{s.t.} \   \|u\|_2=\rho, \ G(u)=0\}$.\\

Let us recall that the equation with a single nonlinearity 
\begin{equation}\label{SP-one}
i\psi_{t}= \sqrt{-\Delta} \pm|\psi|^{p-1}\psi,
\end{equation} 
is  invariant under the scaling $\psi\mapsto \psi_\lambda=\lambda^{\frac{1}{p-1}}\psi(\lambda t, \lambda x)$, and the $L^2$-norm of $\psi_\lambda$ is left invariant provided that $p=1+\frac2d$. Hence, we say that \eqref{SP-one} (and also \eqref{SP}) is mass-subcritical. 
%%%%%%%%%%%%%%%%%%%%%%%%

\subsection{Preliminary results} In this subsection we collect some essential tools which will lead to the proof of the main Theorem \ref{mainthm}.\\

Let us introduce a general scaling $u\mapsto u_\lambda=\lambda^\alpha u(\lambda x)$ often used along the paper. The quadratic part of the energy rescales as
\[
T(u_\lambda)=\lambda^{2\alpha+1-d} T(u),
\]
and a general $L^{r+1}$-norm  rescales as
\[
\|u_\lambda\|_{r+1}=\lambda^{\alpha-d/(r+1)}\|u\|_{r+1}.
\]
Note that by the Plancherel identity, $ \displaystyle T(u)=\int (|\xi|-v\cdot \xi)|\mathcal F(u)|^2d\xi$, then the condition $|v|<1$ ensures that $T(u)>0$ and that $T(u)\sim \|u\|_{\dot H^{1/2}}^2$.  We moreover recall the following Gagliardo-Nirenberg-type inequality: for $1<r<1+\frac{2}{(d-1)_+}$, 
\begin{equation}\label{eq:GN}
\|u\|_{r+1}^{r+1}\leq C(d,v,r) T(u)^{d(r-1)/2}M(u)^{(r+1)/2-d(r-1)/2},
\end{equation}
for any $u\in H^{1/2}$. \\

The first result is the following key Lemma which will play a key role for the rest of the paper.

\begin{lemma}\label{cruciallemma}
Fix $d\in \mathbb N$ and $v\in \R^d$ with $|v|<1$. Let 
\[
A_\delta=\{u \in H^{1/2} \ \hbox{s.t.} \   E(u)\leq 0, \ G(u)=0, \ T(u) \leq \delta \}.
\]
If $\delta$ is sufficiently small then $A_\delta=\emptyset.$
\end{lemma}
\begin{proof}
The first observation  is that in $A_\delta$, the quadratic term $T(u)$ has a size comparable with the focusing term $\|u\|_{p+1}^{p+1}$. Indeed we have
\[
0\geq E(u)-\frac{1}{d(q-1)}G(u)=\frac{d(q-1)-2}{2d(q-1)}T(u)+\frac{p-q}{(p+1)(q-1)}\|u\|_{p+1}^{p+1}
\]
which shows that
\[
T(u)\geq \frac{2d(p-q)}{(2-d(q-1))(p+1)}\|u\|_{p+1}^{p+1}.
\]
On the other hand, the non positivity of the energy implies
that 
\begin{equation}\label{eq-p-control}
\|u\|_{p+1}^{p+1}\geq T(u), 
\end{equation}
therefore 
\begin{equation}\label{eq:norm-equiv}
\|u\|_{p+1}^{p+1}\sim T(u).
\end{equation}
Now, by computing 
$E(u)- \frac{1}{2}G(u)$
we get
\[0\geq E(u)-\frac{1}{2}G(u)=\frac{2-d(q-1)}{2(q+1)}\|u\|_{q+1}^{q+1}-\frac{2-d(p-1)}{2(p+1)}\|u\|_{p+1}^{p+1}.\]
The last inequality guarantees that
\begin{equation}\label{eq:q-p}
\|u\|_{q+1}^{q+1}\lesssim \|u\|_{p+1}^{p+1}.
\end{equation}
At this point, by the Gagliardo-Nirenberg interpolation inequality jointly with the norm equivalence $T(u)\sim \|u\|_{\dot H^{1/2}}^2$, and \eqref{eq:q-p}, we have 
\begin{equation}\label{eq:control-GN}
\begin{aligned}
\|u\|_{p+1}^{p+1}&\leq \|u\|_{q+1}^{(1-\theta)(p+1)}\|u\|_{\dot H^{1/2}}^{\theta(p+1)}\\
&\lesssim \|u\|_{q+1}^{(1-\theta)(p+1)}T(u)^{\theta(p+1)/2}\\
&\lesssim  \|u\|_{p+1}^{(1-\theta)(p+1)^2/(q+1)}T(u)^{\theta(p+1)/2}
\end{aligned}
\end{equation}
with $\theta=\frac{2d(p-q)}{(p+1)(2d-(d-1)(q+1))}$. Combining \eqref{eq-p-control}, \eqref{eq:control-GN}, and \eqref{eq:norm-equiv}, we get
\begin{equation}
T(u)\lesssim  T(u)^{(1-\theta)(p+1)/(q+1)}T(u)^{\theta(p+1)/2}.
\end{equation}
Then, noticing that $1<\frac{(1-\theta)(p+1)}{q+1}+\frac{\theta (p+1)}{2}$ is always verified, as it is equivalent to $2(q-p)<\theta(p+1)(q-1)$, we have that $T(u)$ cannot be too small, namely $A_\delta=\emptyset$ provided that  $
\delta\ll1$.
\end{proof} 

Consider now the energy functional without the defocusing term, i.e., let us introduce 
\[\tilde E(u):=\frac{1}{2}T(u)-\frac{1}{p+1}\|u\|_{p+1}^{p+1},\]
and define
\[
J_{\rho^2}=\inf_{u\in S_{\rho}}\tilde E(u).
\]
Note that when $v=0$, the functional  $\tilde E$ corresponds to the energy functional related to  the equation
\[
i\psi_{t}= \sqrt{-\Delta} \psi -|\psi|^{p-1} \psi,
\]
and the latter is invariant under the scaling $u\mapsto u_\lambda=\lambda^\alpha u(\lambda x)$ with $\alpha=\frac{1}{p-1}$. By this simple observation, we get the following. 
\begin{prop}\label{propomog}
$J_{\rho^2}=\rho^{\frac{4-2(d-1)(p-1)}{2-d(p-1)}}J_1$ and $J_1<0.$
\end{prop}

\begin{proof}
The proof follows by a scaling argument. Let us take $u\in S_1$ and define $u_{\lambda}=\lambda^{\frac{1}{p-1}}u(\lambda x).$  
Straightforward computations give 
\[\|u_{\lambda}\|_{2}^2=\lambda^{\frac{2}{p-1}-d}\|u\|_{2}^2=\lambda^{\frac{2}{p-1}-d},
\] 
thus, by imposing 
$\lambda=\lambda(\rho):=\rho^{\frac{2(p-1)}{2-d(p-1)}}$, we obtain that $u_{\lambda}\in S_{\rho}$.  On the other hand, by the previous observation, $\tilde E(u_{\lambda(\rho)})=\lambda(\rho)^{\frac{2}{p-1}-d+1}\tilde E(u)=\rho^{\frac{4-2(d-1)(p-1)}{2-d(p-1)}}\tilde E(u)$.
The scaling map between  $S_1$ and $S_{\rho}$ is a bijection, hence $J_{\rho^2}=\rho^{\frac{4-2(d-1)(p-1)}{2-d(p-1)}}J_1$.\\
Now we prove that $J_1<0$. Let us consider the mass-preserving scaling $u_{\lambda}=\lambda^{\frac d2}u(\lambda x)$ so that  $u_\lambda \in S_1$ and 
\[
\tilde E(u_{\lambda})= \frac{\lambda}{2}T(u)-\frac {\lambda^{\frac{d(p-1)}{2}}}{p+1}\|u\|_{p+1}^{p+1}.
\]
Note that $1>\frac{d(p-1)}{2}>0$ if and only if $1<p<1+\frac2d$ and that $J_1\leq \tilde E(u_{\lambda})$. The claim follows if we choose $\lambda$ sufficiently small, as $\tilde E<0$ for $\lambda\ll1$.
\end{proof}

We recall the following well-known facts,  that always work for translation invariant minimization problems.
\begin{lemma}\label{lemmainiz}
The function $\rho\rightarrow I_{\rho^2}$ is continuous. Moreover:\\
\textup{(i)}
the ground state energy is weakly subadditive, namely 
\[I_{\rho^2}\leq I_{\mu^2}+I_{\rho^2-\mu^2} \ \text{ for all } \ 0<\mu<\rho;\]
\textup{(ii)} if the ground state energy is strongly  subadditive, namely 
\[I_{\rho^2}<I_{\mu^2}+I_{\rho^2-\mu^2} \ \text{ for all } \  0<\mu<\rho\]
then the infimum is achieved.
\end{lemma}
\begin{proof}
For the proof we refer to the classical reference \cite{Lions-1984}.
\end{proof}

\begin{prop}\label{propasynt}
$J_{\rho^2}\leq I_{\rho^2}\leq0$ and  $\displaystyle \lim_{\rho\rightarrow 0}\frac{I_{\rho^2}}{\rho^{\alpha}}=0$ for all $\alpha \in [2, \frac{4-2(d-1)(p-1)}{2-d(p-1)})$.
\end{prop}
\begin{proof}
The fact that $J_{\rho^2}\leq I_{\rho^2}$ follows from the positivity of the defocusing term. 
The non-positivity of $I_{\rho^2}$ again follows by a scaling argument. Indeed, by taking $u\in S_{\rho}$ and the mass-preserving scaling $u_{\lambda}=\lambda^{\frac d2}u(\lambda x)$ such that $u_{\lambda}\in S_{\rho}$ for all $\lambda>0$, we have
\[E(u_{\lambda})=\frac{\lambda}{2}T(u)+\frac{\lambda^{\frac{d(q-1)}{2}}}{q+1} \|u\|_{q+1}^{q+1}-\frac{\lambda^{\frac{d(p-1)}{2}}}{p+1}\|u\|_{p+1}^{p+1}.
\]
For $\lambda \rightarrow 0$ we have that $E_v(u_{\lambda})\rightarrow 0$, and hence  $ I_{\rho^2}\leq 0$. 
The fact that  that $\displaystyle \lim_{\rho\rightarrow 0}\frac{I_{\rho^2}}{\rho^{\alpha}}=0$ for all $\alpha \in [2,\frac{4-2(d-1)(p-1)}{2-d(p-1)})$ follows from Proposition \ref{propomog} and the fact that $I_{\rho^2}\leq0$.
\end{proof}
We now give a sufficient condition ensuring the strong subadditivity of the function $I_\rho$.
\begin{lemma}\label{lemmonot}
Let $\alpha \in [2, \infty)$ and $I_{\rho^2}<0$. If the function $\displaystyle \rho \rightarrow \frac{I_{\rho^2}}{\rho^{\alpha}}$ is strictly decreasing, then $I_{\rho^2}<I_{\mu^2}+I_{\rho^2-\mu^2}$ for all $\mu \in (0, \rho)$.
\end{lemma}
\begin{proof}
From the inequalities $\frac{\mu^{\alpha}}{\rho^{\alpha}}I_{\rho^2}<I_{\mu^2}$ and $\frac{(\rho^{2}-\mu^{2})^{\frac{\alpha}{2}}}{\rho^{\alpha}}I_{\rho^2}<I_{\rho^2-\mu^2}$ we get, by adding term by term,
\[\frac{\mu^{\alpha}}{\rho^{\alpha}}I_{\rho^2}+\frac{(\rho^{2}-\mu^{2})^{\frac{\alpha}{2}}}{\rho^{\alpha}}I_{\rho^2}<I_{\mu^2}+I_{\rho^2-\mu^2}.\]
As $I_{\rho^2}$ is non-positive by the Proposition \ref{propasynt}, to conclude it suffices to have $\mu^{\alpha}+(\rho^{2}-\mu^{2})^{\frac{\alpha}{2}}<\rho^{\alpha}$. Note that 
\[\mu^{\alpha}+(\rho^{2}-\mu^{2})^{\frac{\alpha}{2}}=(\mu^2)^{\frac{\alpha}{2}}+(\rho^{2}-\mu^{2})^{\frac{\alpha}{2}}<(\rho^2)^{\frac{\alpha}{2}}=\rho^{\alpha}\]
by the convexity inequality $x^a+y^a<(x+y)^a$ if $a\geq 1$.
\end{proof}

The next Lemma shows that a strict negativity of the ground state energy is enough to show the existence of constrained minimizers.
\begin{lemma}\label{negimplcompact}
If $I_{\rho^2}<0$ then exists $u \in S_{\rho}$ such that $E_v(u)=I_{\rho^2}.$
\end{lemma}
\begin{proof}
It suffices from Lemma \ref{lemmonot} that there exists $\alpha\in [2, \infty)$ such that for any $s \in(0, \rho)$,  $\frac{I_{s^2}}{s^{\alpha}}>\frac{I_{\rho^2}}{\rho^{\alpha}}$. Indeed, the monotonicity of $\frac{I_{s^2}}{s^{\alpha}}$ implies strong subadditivity by Lemma \ref{lemmonot}, and the latter implies compactness, see Lemma \ref{lemmainiz}.  See also \cite{BS-JFA}.
\medskip

\noindent Now let us restrict $\alpha $ to the interval $[2, \frac{4-2(d-1)(q-1)}{2-d(q-1)})$, and let us define the quantity 
\[Q_{\alpha}=\inf_{s\in (0, \rho]} \frac{I_{s^2}}{s^{\alpha}}.\]
From the fact that  $I_{\rho^2}<0$, we have that $Q_{\alpha}<0$.  Moreover, Proposition \ref{propasynt} yields 
\[
\rho^{\star}_{\alpha}=\{\inf s \in (0, \rho] \ \hbox{s.t.} \  \frac{I_{s^2}}{s^{\alpha}}=Q_{\alpha}\}>0.
\] 
Clearly, if 
$\rho^{\star}_{\alpha}=\rho$, then the strong  subadditivity holds and a minimizer exists. Therefore, let us assume that
$\rho^{\star}_{\alpha}<\rho.$
In the latter case, we have by definition that  for any $\mu \in (0, \rho^{\star}_{\alpha})$ 
\[
\frac{\mu^{\alpha}}{({\rho^{\star}_{\alpha}})^{\alpha}}I_{({\rho^{\star}_{\alpha}})^2}<I_{\mu^2} 
\]
and 
\[
\frac{((\rho^{\star}_{\alpha})^{2}-\mu^{2})^{\frac{\alpha}{2}}}{\rho^{\alpha}}I_{({\rho^{\star}_{\alpha}})^2}<I_{(\rho^{\star}_{\alpha})^2-\mu^2}.
\]
Hence, by subadditivity, there exists $u_{\alpha}\in S(\rho^{\star}_{\alpha})$ with $E(u_{\alpha})=I_{(\rho^{\star}_{\alpha})^2}$ and such that for $\theta \in (1-\eps, 1+\eps)$, for some small $\eps>0$,
\[
\frac{E(u_{\alpha})}{(\rho^{\star}_{\alpha})^{\alpha}}=\frac{I_{(\rho^{\star}_{\alpha})^2}}{(\rho^{\star}_{\alpha})^{\alpha}}\leq \frac{I_{\theta^2(\rho^{\star}_{\alpha})^2}}{\theta^{\alpha}(\rho^{\star}_{\alpha})^{\alpha}}\leq \frac{E(\theta u_{\alpha})}{\theta^{\alpha}(\rho^{\star}_{\alpha})^{\alpha}}.
\]
Therefore we have 
\begin{equation}\label{eq:condderiv}
 \frac{d}{d\theta} 
\left(\theta^{\alpha} E(u_{\alpha})-E(\theta u_{\alpha})\right)\big\lvert_{\theta=1}=0.
\end{equation}
For a minimizer $u_\alpha$ of $I_{\rho^2}$ we have that 
\begin{equation}\label{eq:rangep}
I_{\rho^2} =E(u_\alpha)=\frac 12T(u_\alpha)+\frac {1}{q+1}\|u_{\alpha}\|_{q+1}^{q+1}-\frac{1}{p+1}\|u_{\alpha}\|_{p+1}^{p+1}
\end{equation}
and 
\begin{equation}\label{eq:rangep-2}
T(u_\alpha)+ \frac{d(q-1)}{q+1} \|u_{\alpha}\|_{q+1}^{q+1}-\frac{d(p-1)}{p+1}\|u_{\alpha}\|_{p+1}^{p+1}=0.
\end{equation}
Observe from \eqref{eq:rangep-2} that 
\begin{equation}\label{eq:c-ab}
\|u_{\alpha}\|_{q+1}^{q+1}=\frac{q+1}{d(q-1)}\left(\frac{d(p-1)}{p+1}\|u_{\alpha}\|_{p+1}^{p+1}-T(u_\alpha)\right),
\end{equation} 
and hence from \eqref{eq:rangep} we get  the identities 
\begin{equation}\label{eq:stimenerg}
I_{\rho^2} =E(u_\alpha)=\frac{d(q-1)-2}{2d(q-1)}T(u_\alpha)+\frac{p-q}{(p+1)(q-1)}\|u_{\alpha}\|_{p+1}^{p+1}
\end{equation}
and
\begin{equation} \label{eq:stimenerg2}
T(u_\alpha)+\|u_{\alpha}\|_{q+1}^{q+1}-\|u_{\alpha}\|_{p+1}^{p+1}=\frac{d(q-1)-q-1}{d(q-1)}T(u_\alpha)+\frac{2(p-q)}{(p+1)(q-1)}\|u_{\alpha}\|_{p+1}^{p+1}.
\end{equation}
From \eqref{eq:condderiv} we obtain
\[\alpha E(u_\alpha)-\left(T(u_\alpha)+\|u_{\alpha}\|_{q+1}^{q+1}-\|u_{\alpha}\|_{p+1}^{p+1}\right)=0,\]
which gives, by using \eqref{eq:c-ab} and  \eqref{eq:stimenerg2},
\begin{equation}\label{eq:a-b}
   T(u_\alpha)=\frac{1}{d(\alpha-2)(q-1)-2q-2-2\alpha}\frac{2d(2-\alpha)(p-q)}{p+1}\|u_{\alpha}\|_{p+1}^{p+1}.
\end{equation}
Note that if $\alpha\in [2, \frac{4-2(d-1)(q-1)}{2-d(q-1)})$, then  $d(\alpha-2)(q-1)-2q-2-2\alpha>0$, while $(2-\alpha)(p-q)$ is always non-positive for $\alpha\geq2$. Hence, for any $\alpha\in[2, \frac{4-2(d-1)(q-1)}{2-d(q-1)})$ we get an absurd from \eqref{eq:a-b} as we would have that $T(u_\alpha)\leq0$.
In the end, $\rho^{\star}_{\alpha}=\rho$ and hence the strong  subadditivity property holds and a minimizer exists.
\end{proof}

In the next Lemma, we show that for a non-empty right neighborhood of $\rho=0$, the infimum of the energy is zero.
\begin{lemma}\label{lem:first} There exists a strictly positive mass $\rho_0$ such that:\\
\textup{(i)} $ I_{\rho^2}=0$  for all  $\rho\in (0,\rho_0]$;\\
\textup{(ii)} $ I_{\rho^2}<0$  for all  $\rho\in (\rho_0, \infty)$.
\end{lemma}
\begin{proof}
The fact that $I_{\rho^2}\leq0$, see Proposition \ref{propasynt}, together with  the weak subadditivity implies that if $I_{\rho^2}<0$ then 
$I_{s^2}<0$ all $s>\rho.$ The  negativity of $I_{\rho^2}<0$ for sufficiently large $\rho$ follows by a scaling argument. Indeed, let us choose  
\[
u_{\lambda}=\lambda^{\frac{1}{q-1}}u(\lambda x);
\]
we have
\[
E(u_{\lambda})=\lambda^{\frac{2}{q-1}-d+1} \left( \frac 12T(u)+\frac{1}{q+1}\|u\|_{q+1}^{q+1}
\right)- \frac{\lambda^{\frac{p+1}{q-1}-d}}{p+1}\|u\|_{p+1}^{p+1}
\]
and then $E(u_{\lambda})<0$ by choosing a  sufficiently large $\lambda$, as $\frac{p+1}{q-1}-d>\frac{2}{q-1}-d+1$ if and only if $p>q$. On the other hand, $\|u_{\lambda}\|_2=\lambda^{\frac{2-d(q-1)}{2(q-1)}}\|u\|_2$ and hence a large $\lambda$  corresponds to a large mass $\rho$. 
 \smallskip
 
\noindent Now we prove (i), i.e., that there exists $\rho_0>0$ such that
$ I_{\rho^2}=0$  for all  $\rho\in ( 0,\rho_0]$. Note that from the weak subadditivity property, together with $I_{\rho^2}\leq 0$, the function $I_{\rho^2}$ is non-increasing. By defining 
\[
\rho_0=\sup\{\rho \ \hbox{s.t.} \   I_{s^2}=0 \text{ for all } s \in (0,\rho)\}
\] 
we have
$I_{\rho^2}<0$ for all $\rho \in (\rho_0, \infty).$
Now we prove that $\rho_0>0.$\\
The idea is to show that $I_{\rho^2}$ cannot be attained in $S_{\rho}$ if $\rho$ is sufficiently small.
As a byproduct we will have that
\begin{equation}\label{rho-0}
\rho_0=\sup\{\rho \ \hbox{s.t.} \  I_{s^2}=0 \text{ for all } s \in (0,\rho]\}
\end{equation}
is strictly positive, because the negativity of $I_{\rho^2}$ implies existence of minimizers thanks to Lemma \ref{negimplcompact}.
\smallskip

\noindent Therefore, let us assume that there exists a sequence 
$\rho_n \rightarrow 0$ such that $I_{\rho_n^2}$ is attained by a ground states $u_{\rho_n}$. By the fact that $E(u_{\rho_n})\leq 0$ we have \eqref{eq-p-control}, so jointly with the Gagliardo-Nirenberg interpolation inequality  we get
\[
T(u_{\rho_n})\lesssim \|u_{\rho_n}\|_{p+1}^{p+1}\lesssim T(u_{\rho_n})^{d(p-1)/2}\rho_n^{p+1-d(p-1)}
\]
and so, as $1<p<1+\frac{2}{d}$, 
\begin{equation}\label{grad-small}
T(u_{\rho_n})=o_n(1).
\end{equation}
On the other hand, a ground state fulfills $G(u_{\rho_n})=0$, hence \eqref{grad-small} contradicts Lemma \ref{cruciallemma}.

\end{proof}

The next Lemma is a non-existence result that shows that if the ground state energy is zero in a open interval then the ground state energy is never achieved. 
\begin{lemma}\label{lemnonex}
If   $I_{\rho^2}=0$ in an interval $\mathcal{I}=(0, \rho_1)$,
then for any $\rho \in \mathcal{I}$, $I_{\rho^2}$ is not attained in $S_{\rho}$.
\end{lemma}
\begin{proof}
Let us assume  that exists $\rho\in \mathcal{I}$
such that $I_{\rho^2}=0=E(u)$ with $u \in S_{\rho}.$
Then 
\[
E(u)=I_{\rho^2} \leq I_{\theta^2\rho^2}\leq E(\theta u)
\]
for $\theta \in (1-\eps,1+\eps)$, for some small $\eps>0$, and then
\[
\frac{d}{d\theta}E(\theta u)\big\lvert_{\theta=1}=0,
\]
which implies
\[
T(u)+\|u\|_{q+1}^{q+1}-\|u\|_{p+1}^{p+1}=0.
\]
The above condition, which tells us  that $u$ is a static solution
solving
\begin{equation}\label{eq:stat}
\sqrt{- \Delta} u +iv\cdot\nabla u   +|u|^{q-1} u-|u|^{p-1}u=0,
\end{equation}
is not compatible with the condition $E(u)=I_{\rho^2}=0.$ Indeed, thanks to \eqref{eq:stimenerg} we get
\[
\|u\|_{p+1}^{p+1}=\frac{(2-d(q-1))(p+1)}{2d(p-q)}T(u)
\]
and thanks to \eqref{eq:stimenerg2}
\[
T(u)+\|u\|_{q+1}^{q+1}-\|u\|_{p+1}^{p+1}=-\frac1dT(u)\neq 0.
\]
This proves that  for any $\rho \in \mathcal{I}$ a minimizer for $E$ on $S_{\rho}$ cannot exist.
\end{proof}

The last Lemma guarantees the existence of a ground state at the critical mass $\rho_0.$
\begin{lemma}\label{lem-mass-cri}
Let $\rho_0$ be defined as in \eqref{rho-0}. Then there exists $u \in S_{\rho_0}$ such that $I_{\rho_0^2}=E(u).$ 
\end{lemma}
\begin{proof}
Let us consider a sequence $\rho_n\rightarrow \rho_0$
with $\rho_n>\rho_0$. We have $I_{\rho_n^2}<0$ and let us call $u_{\rho_n}$ a ground states that belongs to $S_{\rho_n}$. Clearly $u_{\rho_n}$ is a bounded sequence in $H^{1/2}$ and 
$\displaystyle \liminf_{n\to\infty} \|u_{\rho_n}\|_{p+1}^{p+1}>0$.
Indeed, if  along some subsequence $\rho_{n_k}$, $\displaystyle \lim_{k\to\infty}\|u_{\rho_{n_k}}\|_{p+1}^{p+1}=0$,  then by the negativity of the energy  we obtain $\displaystyle \lim_{k\to\infty}T(u_{n_k})=0$, and the latter is in contrast with Lemma \ref{cruciallemma}. By the nonlocal version of the well-known Lieb Translation Lemma,  see \cite{BFV}, up to a space translation and up to a subsequence, $u_{\rho_n}\rightharpoonup \bar u$ with $\bar u \neq 0$, the weak convergence being in $H^{1/2}$. To prove that $\bar u \in S_{\rho_0}$
it suffices to observe that if $\|\bar u\|_2^2=\mu^2<\rho_0^2$ then
\[I_{\rho_0^2-\mu^2}+I_{\mu^2}+o_n(1)\leq E(u_{\rho_n}-\bar u)+E(\bar u)=E(u_{\rho_n})=I_{\rho_0^2}+o_n(1)\]
and hence, by the weak subadditivity inequality, $E(\bar u)=I_{\mu^2}.$ By Lemma \ref{lemnonex} this is a contradiction.
\end{proof}

\subsection{Proof of Theorem \ref{mainthm}} By means of the tools developed above, we can conclude the proof of Theorem \ref{mainthm}. Indeed, it is now an immediate consequence of Lemma \ref{negimplcompact}, Lemma \ref{lem:first},  Lemma \ref{lemnonex}, and Lemma \ref{lem-mass-cri}.

\section{LWP and dynamical results in 1D}

In this section we give a proof of the local well-posedness of \eqref{SP} in the space domain $\R$. Although it follows the same lines of the scheme as in \cite{KLR}, we give a proof here for the sake of clarity and to keep the paper self-contained. 

\subsection{\texorpdfstring{$H^{s}$ solutions, $s>\frac12$}{Hs solutions, s>1/2}} Let us recall the following estimate, see \cite{C-W}, which holds for any function $f\in H^{s}(\R)$, $s>\frac12$:
\begin{equation}\label{eq:Sob}
 \||f|^{p-1}f\|_{H^s}\leq C\|f\|^{p-1}_{\infty}\|f\|_{H^s}.   
\end{equation}
We write a solution to \eqref{SP} in its Duhamel formulation, 
\[
\psi(t)=e^{-it\sqrt{-\Delta}}\psi_0+i\int_0^t e^{-i(t-\tau)\sqrt{-\Delta}}\left( |\psi|^{q-1}\psi- |\psi|^{p-1}\psi \right)(\tau)d\tau,
\]
and we define the ball of radius $R>0$ in the space of functions $C(I,H^s)$, $I=(-T, T)$ being a time interval containing $t=0$, as
\[
B_R=\left\{ f \in L^\infty(I,H^s) \hbox{ s.t. } \|f\|_{L^\infty(I,H^s)}\leq R\right\}.
\]
The latter is a Banach space, and we can perform a fixed point argument in $B_R$ for the map  
\[
\Phi(\psi(t))=e^{-it\sqrt{-\Delta}}\psi_0+i\int_0^t e^{-i(t-\tau)\sqrt{-\Delta}}\left( |\psi|^{q-1}\psi- |\psi|^{p-1}\psi \right)(\tau)d\tau.
\]
By using the unitary property of the linear propagator $e^{-it\sqrt{-\Delta}}$ in any $H^s$, the Minkowski inequality,  the estimate \eqref{eq:Sob}, and the embedding $H^s\hookrightarrow L^\infty$ for $s>\frac12$, we straightforwardly have 
\[
\|\Phi(\psi)\|_{L^\infty(I,H^s)}\leq \|\psi_0\|_{H^s}+CT(\|\psi\|^{q}_{L^\infty(I,H^s)}+\|\psi\|^{p}_{L^\infty(I,H^s)}). 
\]
Hence, for $R=2\|\psi_0\|_{H^s}$ and $T<\frac{1}{2C(R^{q-1}+R^{p-1})}$ we have that $\Phi$ maps $B_R$ into itself. A similar estimate applies for the difference $\Phi(\psi_1(t))-\Phi(\psi_2(t))$, thus $\Phi$ is a contraction in $B_R$ and a solution $\psi(t)$ to \eqref{SP} with $\psi(0)=\psi_0$ exists in $C(I, H^s)$ at least for small times $T>0$. The conservation of energy and mass is a consequence of a standard regularization argument.
The blow-up alternative in $H^{1/2}$ holds also for initial data in $H^s$, as we can see by employing the Brezis-Gallou\"et inequality (see \cite[Appendix D]{KLR}): there exists $C>0$ such that for any $f\in H^s$, $s>\frac12,$
\begin{equation}\label{eq:BG}
\|f\|_{\infty}\leq C\|f\|_{H^{1/2}}\ln^{1/2}\left(2+\frac{\|f\|_{H^s}}{\|f\|_{H^{1/2}}}\right).
\end{equation}
By the fixed point argument above, the blow-up alternative in $H^s$ holds, and precisely  $T_{\max}<\infty$ if and only if $\displaystyle \lim_{t\to T_{\max}^-}\|\psi(t)\|_{H^s}=\infty$, where $T_{\max}$ is the maximal forward time of existence (similarly for the backward time direction). Let u snow introduce the quantity $\displaystyle K=\sup_{t\in I}\|\psi(t)\|_{H^{1/2}}$. By repeating the estimate of the contraction argument, and by exploiting \eqref{eq:BG}, we have that 
\[
\|\psi(t)\|_{H^s}\leq \|\psi_0\|_{H^s}+C\int_0^t\|\psi(\tau)\|_{H^{1/2}}^{p-1} \ln^{{\frac{p-1}{2}}}\left( 2+\frac{\|\psi(\tau)\|_{H^s}}{\|\psi(\tau)\|_{H^{1/2}}}\right)\|\psi(\tau)\|_{H^s}d\tau,
\]
and by noting that $g(x)=x^{p-1}\ln^{\frac{p-1}{2}}\left( 2+a/x\right)$ is monotone increasing, we get  
\[
\|\psi(t)\|_{H^s}\leq \|\psi_0\|_{H^s}+CK^{p-1}\int_0^t\ln^{{\frac{p-1}{2}}}\left( 2+\frac{\|\psi(\tau)\|_{H^s}}{K}\right)\|\psi(\tau)\|_{H^s}d\tau.
\]
By setting $h(t)=\frac{\|\psi(t)\|_{H^s}}{K}$,  we rewrite the estimate above as  
\[
h(t)\leq \frac{\|\psi_0\|_{H^s}}{K}+CK^{p-1}\int_0^t\ln^{{\frac{p-1}{2}}}( 2+h(\tau))h(\tau)d\tau:=H(t).
\]
Observe that 
\[
\frac{d}{dt}H(t)=K^{p-1}h(t)\ln^{{\frac{p-1}{2}}}( 2+h(t))\leq K^{p-1}(2+H(t))\ln^{{\frac{p-1}{2}}}(2+H(t)),
\]
and then
\[
\frac{d}{dt}\left(\ln^{\frac{3-p}{2}}(2+H(t))\right)\leq \frac{3-p}{2}K^{p-2} 
\]
which easily gives 
\[
2+H(t)\leq \exp \left(\frac{3-p}{2}K^{p-2}t+\ln^{\frac{3-p}{2}}\left(2+\frac{\|\psi_0\|_{H^s}}{K}\right) \right)^{\frac{2}{3-p}}
\]
and then the norm $\|\psi(t)\|_{H^s}$ remains bounded if $T$ is bounded. This means that the blow-up alternative holds at the regularity $H^{1/2}$.
At this point, by using the blow-up alternative in $H^{1/2}$, and the mass-subcritical nature of the nonlinearities, we can infer that the solutions can be extended globally in time.
%({\color{red} CONTINUE})

\subsection{\texorpdfstring{$H^{1/2}$ solutions}{H1/2 solutions}}
The existence of solutions at the $H^{1/2}$ regularity  uses a regularization argument, see \cite{KLR} (see also  \cite{Chen} for a Half-Wave-Schr\"odinger equation in $\R^2$). Let us consider a sequence of initial data $\{\psi_{0,n}\}\subset H^{s}$, $s>\frac12$, converging to $\psi_0$ in the $H^{1/2}$-topology. The $H^s$ solutions $\psi_n(t)$ to \eqref{SP} corresponding to these initial data are global by the previous discussion, and moreover the $H^{1/2}$-norm of these solutions remains bounded by the conservation of the energy, as well as remains bounded the $H^{-1/2}$-norm of $\partial_t\psi_n(t)$. The bound is actually uniform in $n$, and then locally in time we can extract a weakly convergent subsequence $\psi_n(t)\rightharpoonup \psi(t)$, and by compact embedding the convergence is actually strong in the $L^r_{\loc}$-topology and then $\psi$ is a weak solution to \eqref{SP}.

By employing an argument to due to Judovic \cite{Ju}, see  \cite{GG} and \cite{KLR} for the application to the Szeg\H{o} equation and to the cubic Half-Wave equation, respectively,   we prove uniqueness of weak solutions. Let   $\psi_1$ and $\psi_2$ be two solutions emanating from the same initial datum $\psi_0$, and define the function $h(t)=\|\psi_1(t)-\psi_2(t)\|_2^2$. Using the equation solved by $\psi_1$ and $\psi_2$, by the fact that $\Im \int( \psi_1- \psi_2) \sqrt{-\Delta}( \psi_1- \psi_2)=0$, we have 
\begin{equation}
    \begin{aligned}
        \frac{d}{dt}h(t)&=2\Re \int(\partial_t \psi_1-\partial_t \psi_2)(\bar \psi_1-\bar \psi_2)\\
        &=2\Im \int (\bar \psi_1-\bar \psi_2)(|\psi_2|^{p-1}\psi_2-|\psi_1|^{p-1}\psi_1)\\
        &+ 2\Im \int (\bar \psi_1-\bar \psi_2)(|\psi_2|^{q-1}\psi_2-|\psi_1|^{q-1}\psi_1).
    \end{aligned}
\end{equation}
The analysis of the term involving $(|\psi_2|^{\bullet-1}\psi_2-|\psi_1|^{\bullet-1}\psi_1)$, $\bullet\in\{p,q\}$ is analogous, so we consider just one of them. By the following easy computations, we have that 
\begin{equation}\label{eq:usefull} 
\begin{aligned}
&\phantom{=}|\psi_1-\psi_2|^2(|\psi_1|^{p-1}+|\psi_2|^{p-1})=|\psi_1-\psi_2|^{2-\delta}|\psi_1-\psi_2|^{\delta}(|\psi_1|^{p-1}+|\psi_2|^{p-1})\\
&\lesssim |\psi_1-\psi_2|^{2-\delta}(|\psi_1|^\delta +|\psi_2|^\delta)(|\psi_1|^{p-1}+|\psi_2|^{p-1}) \\
&\lesssim   |\psi_1-\psi_2|^{2-\delta} (|\psi_1|^{p-1+\delta}+|\psi_2|^{p-1+\delta}+|\psi_2|^{\delta}|\psi_1|^{p-1}+|\psi_1|^{\delta}|\psi_2|^{p-1})\\
&\lesssim |\psi_1-\psi_2|^{2-\delta} (|\psi_1|^{p-1+\delta}+|\psi_2|^{p-1+\delta}),
\end{aligned}
\end{equation}
for any $\delta\in(0,2)$, where in the last step we used the Young inequality with exponents $p=\frac{p-1+\delta}{\delta}$ and $p'=\frac{p-1+\delta}{p-1}$. By using in order $
||\psi_1|^{p-1}-|\psi_2|^{p-1}|\lesssim |\psi_1-\psi_2|(|\psi_1|^{p-1}+|\psi_2|^{p-1})$,
\eqref{eq:usefull}, and the H\"older inequality, we get  
\begin{equation}
    \begin{aligned}
        &2\Im \int (\bar \psi_1-\bar \psi_2)(|\psi_2|^{p-1}\psi_2-|\psi_1|^{p-1}\psi_1)\\
        &\lesssim \int |\psi_1-\psi_2|^2(|\psi_1|^{p-1}+|\psi_2|^{p-1})\\
        &\lesssim \int |\psi_1-\psi_2|^{2-\delta}(|\psi_1|^{p-1+\delta}+|\psi_2|^{p-1+\delta})\\
        &\lesssim \|\psi_1-\psi_2\|_{L^{r(2-\delta)}}^{2-\delta}(\|\psi_1\|_{L^{r'(p-1+\delta)}}^{p-1+\delta} +\|\psi_2\|_{L^{r'(p-1+\delta)}}^{p-1+\delta}).
    \end{aligned}
\end{equation}
Choosing $r=\frac{2}{2-\delta}$ we then have
\begin{equation}
    \begin{aligned}
        &2\Im \int (\bar \psi_1-\bar \psi_2)(-|\psi_1|^{p-1}\psi_1+|\psi_2|^{p-1}\psi_2)\\
        &\lesssim \|\psi_1-\psi_2\|_{L^{2}}^{2/r}(\|\psi_1\|_{L^{r'(p-1+\delta)}}^{p-1+\delta} +\|\psi_2\|_{L^{r'(p-1+\delta)}}^{p-1+\delta}).
    \end{aligned}
\end{equation}
We recall now that there exists a constant $C>0$ such that for any $1<r<\infty$, $\|f\|_r\leq C\sqrt r\|f\|_{H^{1/2}}$, for any $f\in H^{1/2}$, see \cite[Appendix D]{KLR}. It follows that 
\begin{equation}
\begin{aligned}
        \frac{d}{dt}h(t)&\lesssim (r'(p-1+\delta))^{\frac{p-1+\delta}{2}}h(t)^{1/r}(\|\psi_1\|_{H^s}^{p-1+\delta} +\|\psi_2\|_{H^s}^{p-1+\delta})\\
        &\lesssim(r'(p-1+\delta))^{\frac{p-1+\delta}{2}}h(t)^{1/r}
        \end{aligned}
\end{equation}
where in the last step we used the uniform bound $\displaystyle\sup_t(\|\psi_1(t)\|_{H^s}+\|\psi_2(t)\|_{H^s})\leq C$, which implies 
\[
\frac{d}{dt}\left(r'h(t)^{1/r'} \right)\lesssim (r'(p-1+\delta))^{\frac{p-1+\delta}{2}}
\]
and then, for any $t$,
\[
h(t)\lesssim \left(\frac{t}{r'}(r'(p-1+\delta))^{\frac{1}{2}} \right)^{(p-1+\delta)r'}=\left(\frac{t\sqrt{p-1+\delta}}{\sqrt{r'}}\right)^{(p-1+\delta)r'}\to 0
\]
as $r'\to \infty$. Note that $r'\to\infty$ if and only if $\delta\to0$. Hence, the uniqueness of weak solutions is proved. 

We can pass from weak to strong solutions by the following argument. We know that 
\[
\|\psi(t)\|_{2}\leq \liminf_{n\to\infty}\|\psi_n(t)\|_2=\liminf_{n\to\infty}\|\psi_{0,n}\|_2=\|\psi_0\|_2,
\]
where we used the weak convergence of $\psi_n(t)$ to $\psi(t)$, the conservation of the mass for the $\psi_n(t)$ and the strong convergence of $\psi_{0,n}$ to $\psi_0$, in $L^2$. Since the non-linear flow \eqref{SP} is time reversible, one gets the converse inequality, and thus $\|\psi(t)\|_2=\|\psi_0\|_2$. By this fact,
\[
\|\psi(t)\|_2=\lim_{n\to\infty}\|\psi_{0,n}\|_2=\lim_{n\to\infty}\|\psi_{0,n}(t)\|_2,
\]
always by conservation of the mass, and then $\lim_{n\to\infty}\|\psi_{n}(t)\|_2=\|\psi(t)\|_2$.
This implies that locally in time, we have strong convergence $\psi_n(t)\to\psi(t)$ in the $L^2$-topology. By invoking the Gagliardo-Nirenberg inequality, and the uniform bound of the $H^{1/2}$-norm of the weak solutions, uniformly in time, we obtain that $\psi_n(t)\to\psi(t)$ strongly in $L^{q+1}\cap L^{p+1}$ as well, locally in time. By strong convergence in $ L^{q+1}\cap L^{p+1}$ and the conservation of mass and energy, we have
\[
\begin{aligned}
&\frac{1}{2}\|\psi(t)\|_{H^{1/2}}^2+\frac{1}{q+1}\|\psi(t)\|_{q+1}^{q+1}-\frac{1}{p+1}\|\psi(t)\|_{p+1}^{p+1}\\
&\leq\lim_{n\to\infty}\left(\frac{1}{2}\|\psi_n(t)\|_{H^{1/2}}^2+\frac{1}{q+1}\|\psi_n(t)\|_{q+1}^{q+1}-\frac{1}{p+1}\|\psi_n(t)\|_{p+1}^{p+1}\right)\\
&=\lim_{n\to\infty}\left(\frac{1}{2}\|\psi_{0,n}\|_{H^{1/2}}^2+\frac{1}{q+1}\|\psi_{0,n}\|_{q+1}^{q+1}-\frac{1}{p+1}\|\psi_{0,n}\|_{p+1}^{p+1}\right)\\
&=\frac{1}{2}\|\psi_{0}\|_{H^{1/2}}^2+\frac{1}{q+1}\|\psi_{0}\|_{q+1}^{q+1}-\frac{1}{p+1}\|\psi_{0}\|_{p+1}^{p+1}.
\end{aligned}
\]
By the time reversibility of the equation and uniqueness, again, we find the converse inequality for any $t$, and then, locally uniformly in time, we have $\lim_{n\to\infty}\|\psi_n(t)\|_{H^{1/2}}=\lim_{n\to\infty}\|\psi(t)\|_{H^{1/2}}$, then we upgrade the weak converge of $\psi_n(t)\rightharpoonup \psi(t)$ to a strong convergence, i.e., $\psi_n(t)\to\psi(t)$ in $H^{1/2}$. This completes the 
existence theory in $H^{1/2}$.

\subsection{Stability result} We conclude with the proof of Corollary \ref{thm:corol}.
    In this section, we prove Corollary \ref{thm:corol} following the ideas of \cite{CL-CMP}. Fix $d=1$,  $v\in \R$ with $|v|<1$, and  let $\rho_0^v$ be as given in Theorem \ref{mainthm}. For any $\rho>\rho_0^v$, the set 
\[
\mathcal G=\{e^{i\gamma}u(\cdot+y) \ \hbox{s.t.} \   \gamma\in\R,  \ y\in \R, \  \hbox{ and } \ u\in S_\rho \ \hbox{ with } \ I_{\rho^2}=E_v(u)\}
\]is \emph{orbitally stable} if
for every $\varepsilon>0$ there exists $\delta>0$ such that for any $\psi_{0}\in H^{\frac 12}(\R)$ with
$\inf_{w\in \mathcal G}\|w-\psi_{0}\|_{H^{\frac 12}}<\delta$ we have
\[
\sup_{t\in \R}\inf_{w\in \mathcal G
} \|\psi(t,\cdot)-w\|_{H^{1/2}(\R)}<\varepsilon
\]
where $\psi(t,\cdot)$ is the solution of \eqref{SP} with initial datum $\psi_{0}$. Note that $\mathcal G$ is invariant by translation, namely 
if $w\in \mathcal G$ then also $w(\cdot-y)\in \mathcal G$ for any $y\in \R.$\\
We also notice that  our functional $E_{v}$ is invariant along the time evolution.
We argue by contradiction
assuming that there exists a $\rho>\rho_{0}^v$ such that $\mathcal G$ is not orbitally stable. This means that there exists $\varepsilon>0$, a sequence of initial
data $\{\psi_{n,0}\}\subset H^{1/2}(\R)$, and a sequence of times  $\{t_{n}\}\subset\R$, such that the global  solution $\psi_{n}(t)=$ with
$\psi_{n}(0,\cdot)=\psi_{n,0}$, satisfies
\begin{equation*}
\lim_{n\rightarrow \infty}\inf_{w\in \mathcal G}\|\psi_{n,0}-w\|_{H^{1/2}(\R)}=0 \quad \text{ and }\quad \inf_{w\in \mathcal G}\|\psi_{n}(t_{n},\cdot)-w\|_{H^{1/2}(\R)}\ge\varepsilon.
\end{equation*}
Then there exists a  minimizer $u_{\rho}\in H^{1/2}(\R)$ of $E_v$ and $\theta\in \R$ such that $w=e^{i\theta}u_{\rho}$
and
\[
\|\psi_{n,0}\|_{2}\rightarrow\|v\|_{2}=\rho \quad \text{ and } \quad
E_v(\psi_{n,0})\rightarrow E_v(u_{\rho}).
\]
Note that we can assume that
$\psi_{n,0}\in S_{\rho}$. Indeed, by setting $\alpha_{n}=\rho/\|\psi_{n,0}\|_{2}$, we have that $\alpha_{n}\rightarrow1$, 
$\alpha_{n}\psi_{n,0}\in S_{\rho}$ and $ E_v(\alpha_{n}\psi_{n,0}) 
\rightarrow I_{\rho^2}$, and then we can replace $\psi_{n,0}$ with
 $\alpha_{n}\psi_{n,0}$.  Hence,   $\{\psi_{n,0}\}$ is a minimizing sequence for $I_{\rho^2}$, and since
 $E_v(\psi_{n}(t_{n}))=E_v(\psi_{n,0}),$
also $\{\psi_{n}(t_{n}, \cdot)\}$ is a minimizing sequence for $I_{\rho^2}$. Since we proved that every minimizing sequence has a converging subsequence  (up to translation) in the strong $H^{\frac 12}$-topology  to a minimum on the sphere $S_{\rho}$, 
we have a contradiction and the proof is complete.

\subsection*{Acknowledgements}\rm 
J. B. acknowledges  Louis  Jeanjean for useful discussions on the mass-subcritical NLS equation, and for pointing-out reference \cite{JL-2022}. L. F. is member of the GNAMPA of the INdAM (Instituto Nazionale di Alta Matematica).

\end{document}